\documentclass[a4paper,11pt]{article}
\usepackage{geometry,latexsym,amssymb,amsmath,color}
\usepackage[latin5]{inputenc}
\usepackage{enumerate}
\usepackage{enumitem}
\usepackage{hyperref}
\usepackage[all]{xy}
\usepackage{authblk}
\usepackage{tikz}
\usepackage{tikz-cd}
\usepackage{amsthm}
\usepackage{amsfonts}

\usetikzlibrary{matrix,arrows}

\newtheorem{example}{Example}[section]
\newtheorem{Def}[example]{Definition}
\newtheorem{Exam}[example]{Example}
\newtheorem{Prop}[example]{Proposition}
\newtheorem{Thm}[example]{Theorem}
\newtheorem{Lem}[example]{Lemma}
\newtheorem{Rem}[example]{Remark}

\newenvironment{Prf}{{\bf Proof:} }{\hfill $\Box$\mbox{}}

\def\Cov{\mathsf{Cov}}

\def\Cat{\mathsf{Gpd}}
\def\Act{\mathsf{Act}}

\def\St{\mathsf{St}}
\def\C{\mathsf{Set^\mathbb{T}}}

\def\leq{\leqslant}

\def\P{\mathsf{P}}

   \def\Top{\mathsf{Top}}
\def\T{\mathbb{T}}
\def\Set{\mathsf{Set}}
\def\Ob{\mathsf{Ob}}

\def\Alg {\mathsf{\Alg}}

\def\leq{\leqslant}

\begin{document}
\title{\large\bf Covering morphisms of internal groupoids in the models of a  semi-abelian theory}
\author[1]{Osman MUCUK \thanks{mucuk@erciyes.edu.tr}}
\author[2]{Serap DEM\.{I}R \thanks{srpdmr89@gmail.com}}
\affil[1]{Department of Mathematics, Erciyes University Kayseri 38039, Turkey}
\affil[2]{Department of Mathematics, Erciyes University Kayseri 38039, Turkey}
\date{\vspace{-5ex}}
\maketitle

\\ {\bf Classification:} 18D35, 18B30, 22A05, 57M10



\begin{abstract} In this paper,  for given an algebraic theory $\T$ whose category $\C$ of models is semi-abelian,  we consider the topological models of $\T$ called topological $\T$-algebras and obtain some results related to the fundamental groups of  topological $\T$-algebras.  We also deal with the internal groupoid structure in the category of models providing that the fundamental groupoid deduces  a functor from topological $\T$-algebras to the internal groupoids in $\C$  and  prove a criterion  for the lifting of such an internal groupoid structure to the covering groupoids.
\end{abstract}

\section{Introduction}

Semi-abelian categories are Barr-exact, protomodular which means the short five lemma holds, have finite coproducts and have zero objects  \cite{Janelidze} . For example all abelian categories, the category of all groups, of rings without unit, of $\Omega$-groups, of Heyting semi-lattices, of locally boolean distributive lattices, of loops, of presheaves or sheaves of these are semi-abelian categories.

Let $\T$ be an algebraic theory in the sense of   Lawvere \cite{Law}. A model of $\T$ in the category of sets is called {\em $\T$-algebra} and a model of $\T$ in the category of topological spaces is called {\em topological $\T$-algebra}.  An algebraic theory $\T$ whose category $\Set^{\T}$ of the  models is semi-abelian is called {\em semi-abelian theory}  and  a model of such a theory is called {\em semi-abelian algebra}.  Such a theory is characterized in  \cite[Theorem 1.1]{Brn2}.  The category of topological $\T$-algebras and continuous $\T$-homomorphisms between them is denoted by $\Top^{\T}$. For example  when  $\T$ is the theory of groups,  $\Top^{\T}$ becomes the category of topological groups.

In  this paper for a semi-abelian theory $\T$ we prove that  the fundamental group of a topological $\T$- algebra is a $\T$-algebra and obtain topological $\T$-algebras  corresponding to the subalgebras.    We also prove that the fundamental groupoid of a topological $\T$-algebra is an internal groupoid in $\C$ and obtain some results on covering morphisms of internal groupoids.

On the one hand in \cite{Bor-Man} some properties of topological groups such as being Hausdorff, compact, connected and etc. have been generalized to the topological $\T$-algebras for a semi-abelian algebra.  Likewise  similar  results have been obtained in \cite{Bor-Man2} for more wider class of algebras called as  topological protomodular algebras which have $n$-constants rather than unique constant.

On the other hand we know from \cite[Theorem 10.34]{Rot} that for  a  connected topological space $X$ which has a universal cover, $x_0\in X$ and a subgroup $G$ of the fundamental group $\pi_1(X,x_0)$ at $x_0$, there is  a covering map $p\colon
(\widetilde{X}_G,\tilde{x}_0)\rightarrow (X,x_0)$ of pointed spaces, with characteristic group $G$ and by \cite[Theorem 10.42]{Rot} whenever  $X$ is a topological group,  $\widetilde{X}_G$ becomes a topological group such that $p$ is a group homomorphism. Using this method  Mucuk and \c{S}ahan in \cite{Mu-Sa} recently have  generalized some results on covering groups of topological groups to the topological groups with operations whose idea comes from Higgins \cite{Hig2},  Orzech  \cite{Orz1} and   \cite{Orz2}.

For a given semi-abelian theory $\T$ and a topological $\T$-algebra $A$, we start our work by proving that the fundamental group $\pi_1(A,e)$ at the constant  $e\in A$ becomes a $\T$-algebra. Then assuming that  $A'$  is a sub $\T$-algebra  of $A$ and $B$ is the fundamental group  $\pi_1(A',e)$, we prove that $\widetilde{A}_B$  corresponding to $B$ as in  \cite[Theorem 10.34]{Rot} is a topological $\T$- algebra such that  $p\colon (\widetilde{A}_B,\tilde{e})\rightarrow (A,e)$ is a topological $\T$-homomorphism.
Next we define internal category in the category $\Set^{\T}$ of semi-abelian algebras and prove that the fundamental groupoid $\pi A$ of a topological $\T$-algebra $A$ is an internal groupoid in   $\Set^{\T}$.  We continue  proving an  equivalence of the categories in Theorem \ref{Covandactareequval};  and finally we obtain  a criterion   for the lifting of internal groupoid structure to the covering groupoids considering the internal groupoid structure in the category  $\C$.

We acknowledge that an extended abstract including  the statements without proofs of some results of this paper appears in \cite{Mu-De} as AIP
proceedings of the conference.

\section{Preliminaries on  covering groupoids and  topological semi-abelian algebras}\label{PrelimCov}

A {\em groupoid}  is a small category in which each arrow  is an isomorphism  (see  \cite{Br} and \cite{Ma} for more discussion on groupoids). More precisely   a groupoid $G$ has a set $G$ of arrows   and  a set $\Ob(G)$ of {\em objects} together with {\em source} and {\em target} point maps $s, t\colon G\rightarrow \Ob(G)$  and {\em object inclusion} map  $\epsilon \colon \Ob(G) \rightarrow G$  such that $s\epsilon=t\epsilon=1_{\Ob(G)}$. There exists a partial composition defined by $ G_t\times_s G\rightarrow G, (g,h)\mapsto g\circ h$, where $G_t\times_s G$ is the pullback of $t$ and $s$.  Here if $g,h\in G$ and $t(g)=s(h)$, then the {\em composite}  $g\circ h$ exists such that $s(g\circ h)=s(g)$ and $t(g\circ h)=t(h)$. Further, this partial composition is associative, for $x\in \Ob(G)$ the arrow  $\epsilon (x)$ acts as the identity and it is denoted by $1_x$, and each arrow $g$ has an inverse $g^{-1}$ such that $s(g^{-1})=t(g)$,  $t(g^{-1})=s(g)$,   $g\circ g^{-1}=\epsilon (s(g))$, $g^{-1}\circ g=\epsilon (t((g))$. The map $G\rightarrow G$, $g\mapsto g^{-1}$ is called the {\em inversion}. In  a groupoid $G$, the source and target points, the object inclusion, the inversion maps and the partial composition  are  called {\em structural maps}.

An example of a groupoid is the {\em fundamental groupoid} $\pi(X)$  of a topological space $X$, where the objects are points of $X$ and  arrows, say from $x$ to $y$ are the  homotopy classes of the paths  in $A$, relative to the end points, with source point $x$ and final point $y$. The partial composition on the homotopy classes is defined by the concatenation of the  paths.  A group is also a groupoid with one object.

In a groupoid $G$ for $x,y\in \Ob(G)$ we write $G(x,y)$ for the set of all
arrows with source points $x$ and target points $y$. According to \cite{Br} $G$ is
{\em transitive} if for all $x,y \in \Ob(G)$, the set  $G(x,y)$ is not empty; for
$x \in \Ob(G)$  the {\em star}  of $x$ is defined  as $\{g\in G\mid s(g)=x\} $ and denoted by  $\St_G x$; and the {\em object group} $G(x,x)$  at $x$  is  denoted by $G(x)$.

A functor  $p\colon H\rightarrow G$ of groupoids is called a {\em morphism} of groupoids. A groupoid morphism $p\colon H\rightarrow G$ is said to be {\em covering morphism}  and $H$   {\em  covering groupoid} of $G$ if for
each $\tilde x\in \Ob(H)$ the restriction  $\St_{H} {\tilde{x}}\rightarrow \St_G{p(\tilde x)}$  is  bijective.   A covering morphism $p\colon H\rightarrow G$ in which  both $H$ and $G$ are transitive is called {\em universal} when $H$ covers every cover of $G$ in the sense that  for every covering morphism $q\colon K\rightarrow G$ there is a unique  morphism of groupoids $r\colon H\rightarrow K$ such that $r\circ q=p$.

For a groupoid morphism $p\colon H\rightarrow G$ and
an object $\tilde{x}$ of $H$ we call  the subgroup
$p(H(\tilde{x}))$ of $G(p\tilde{x})$  as  {\em
characteristic group} of $p$ at $\tilde{x}$.

An {\em action} of a groupoid  $G$ on a set $A$ is defined  in \cite[pp.373]{Br} as consisting of two functions $\omega\colon A\rightarrow \Ob(G)$ and $\varphi\colon A_{\omega}\times_{s}G\rightarrow A,
(a,g)\mapsto ag$, where  $A_{\omega}\times_{s} G$ is the pullback of $\omega$ and  $s$,  subject to the following conditions:
\begin{enumerate}
\item $\omega(ag)=t(g)$ for $(a,g)\in A_{\omega}\times_{s} G$;

\item  $a(g\circ h)=(ag)h$ for  $(g,h)\in G_{t}\times_{s}G$ and $(a,g)\in A_{\omega}\times_{s} G $;

\item  $a\epsilon({\omega(a))}=a$ for $a\in A$.
\end{enumerate}

According to \cite[10.4.2]{Br}  for given such an action of groupoid  $G$  on a set $A$, the {\em semi-direct product groupoid } $G\ltimes A$ with  object set $A$ is defined such that the arrows from  $a$ to $b$ are  the pairs $(g,a)$ with $g\in
G(\omega(a),\omega(b))$ and $ag=b$. The partial composition is given  by
\[  (g,a)\circ (h,b)=(g\circ h,a)\] when $b=ag$. The projection map $p\colon G\ltimes A\rightarrow G$ defined on objects by $\omega$   and on arrows by $(g,a)\mapsto g$ is a covering morphism.

Since in the proof of Theorem \ref{semiabelactiongpd}, we need some details of the   following result
we remind  a  sketch proof from \cite[10.4.3]{Br}.

\begin{Thm}  \label{Theoactiongpdcover} Let $x$ be an object of a transitive  groupoid $G$,
and let $C$ be a subgroup of the object group $G(x)$. Then there exists a covering morphism
$q\colon (H_C,\tilde{x})\rightarrow (G,x)$ with
characteristic group $C$. \end{Thm}
\begin{Prf} Let $A_C$ be the set of  cosets
$C\circ g=\{c\circ g\mid c\in\ C\}$ for $g$ in $St_G x$. Let $\omega\colon
A_C\rightarrow \Ob(G)$ be a map, which maps  $C\circ g$ to the target point of $g$. The function $\omega$ is well defined, because  if $C\circ g=C\circ h$ then  $t(g)=t(h)$.  The groupoid $G$ acts on $A=A_C$ by \[\varphi\colon A_{\omega}\times_{s}G\rightarrow A, (C \circ g,
h)\mapsto C\circ (g\circ h).\] The required groupoid $H_C$ is
taken to be the semi-direct product groupoid $G\ltimes A_C$. Then the  projection
$q\colon H_C\rightarrow G$ given on objects by
$\omega\colon A_C\rightarrow \Ob(G)$ and on arrows by $(h,C \circ g)\mapsto
h$, is a covering morphism of groupoids and has the characteristic
group $C$.  Here  the partial composition  on $H_C$ is defined by
\[(k,C\circ g)\circ(l,C\circ h)=(k\circ l,C\circ g)\]
whenever $C\circ h=C\circ g \circ k $. The required object
$\tilde{x}\in H_C$ is the coset $C$.
\end{Prf}

An {\em algebraic theory} $\T$ in the  sense of  Lawvere \cite[pp.109]{Law}  is a category with objects \[T^0,T, T^2, T^3,\dots \]  where $T^n$ is $n$-copies of the {\em distinguish object}  $T$, and with $m$ arrows
\[\pi_i ^{(m)}\colon T^{m}\rightarrow T, i=0,1,2,\dots,m-1\] for each $m$ such that for  any $m$ arrows \[\tau_i\colon T^{n}\rightarrow T, i=0,1,\dots,m-1 \] in $\T$ there is exactly one arrow \[(\tau_0,\tau_1,\dots,\tau_{m-1})\colon T^{n}\rightarrow T^{m}\] so that
\begin{align*}\label{TalgebraAxiom}
(\tau_0,\tau_1,\dots,\tau_{m-1})\circ \pi_i ^{(m)}=\tau_i ~~ (i=0,1,\dots,m-1)
\end{align*}

The arrows $\tau\colon T^{n} \rightarrow T$ of this category are called {\em $n$-ary operations} and in particular, the $0$-ary operations $T^0\rightarrow T$ are called {\em constants} of the theory  $\T$ .

Throughout the paper by an algebraic theory $\T$ we mean the theory in the sense of Lawvere as stated (see also   Bourceux \cite{Borceux} for an equivalent  set theoretical interpretation).

A product preserving functor $F\colon \T\rightarrow \Set$ is called a {\em model} of the theory or a {\em $\T$-algebra} and  natural transformations between $\T$-algebras are called {\em $\T$-homomorphisms}. Hence a $\T$-homomorphism is a map between sets commuting with all operations of the theory. Let  $\C$ be the  category of models of the algebraic theory $\T$ whose objects are  $\T$-algebras and arrows are  $\T$-homomorphisms.

   A semi-abelian theory  is characterized in  \cite[Theorem 1.1]{Brn2} as follows.  We write $e$  rather than $0$ for the constant of the theory $\T$ to  distinguish from $0$ in a path $\beta \colon [0,1]\rightarrow A$ for a topological $\T$-algebra $A$.
   \begin{Thm} \label{Theosat}
An algebraic theory $\mathbb{T}$ has a semi-abelian category $\C$ of models precisely when,  for some natural number $n$,  the theory $\mathbb{T}$ contains
\begin{enumerate}
\item a unique constant $e$;
\item $n$ binary operations $\alpha_1 (X,Y)$, $\alpha_2 (X,Y)$,\dots,$\alpha_n (X,Y)$ satisfying  $\alpha_i (X,X)=e$;
\item \label{sattheta} an $(n+1)$-ary operation $\theta(X_1,X_2,\dots,X_{n+1})$ satisfying \[\theta(\alpha_1 (X,Y),\alpha_2 (X,Y),\dots,\alpha_n (X,Y),Y)=X.\]
\end{enumerate}
\end{Thm}

Here we  remark  that, in general, $\mathbb{T}$ admits many more operations than simply $\alpha_i$ and $\theta$;  and the choice in $\mathbb{T}$ of the operations  $\alpha_i$ and $\theta$ as indicated is not unique. We mean such a theory by {\em semi-abelian theory} and the corresponding    $\T$-algebras by   {\em semi-abelian algebras}.

For example; each algebraic theory $\mathbb{T}$ which has a unique constant and a group operation  `$+$' is semi-abelian. This is in particular the case for groups, abelian groups,  $\Omega$-groups, modules on a ring,  rings or algebras without unit, Lie algebras,   etc. In Theorem \ref{Theosat} one  chooses $n=1$ and  $\alpha_1 (X,Y)=X-Y$, $\theta(X,Y)=X+Y$.

Replacing $\Set$ with $\Top$ in the definition of $\T$-algebra,  we obtain the categorical  notion of topological $\T$-algebra as a functor
$F\colon \T\rightarrow \Top$.  An equivalent set theoretical definition of a topological $\T$-algebra is given in \cite[Definition 5]{Bor-Man2} as follows.

\begin{Def}\label{TopologicalTalgebra} {\em   Given an algebraic theory $\mathbb{T}$, by a {\em topological $\mathbb{T}$-algebra} we mean a topological
space $A$ provided with the structure of a $\mathbb{T}$-algebra, in such a way that every operation
$\tau\colon T^n\rightarrow T$ of $\T$ induces a continuous mapping
$\tau_A\colon A^n\rightarrow A, (a_1, \dots, a_n)\mapsto \tau(a_1, \dots, a_n)$.
}\end{Def}

We write $\Top^{\mathbb{T}}$ for the category of topological $\mathbb{T}$-algebras and continuous $\mathbb{T}$-homomorphisms
between them.
If $\mathbb{T}$ is a semi-abelian theory, the corresponding topological $\mathbb{T}$-algebras will be also called {\em topological  semi-abelian  algebras}.

We recall from \cite{Bor-Man2} that a theory $\T$ is called {\em protomodular}  if the category $\C$ of models of the theory is protomodular as defined by Bourn in \cite{BournNormalization} and a protomodular theory generalizing semi-abelian theory  is characterized by Bourn and Janelidze in \cite{Brn2} as a theory with $n$-constants $e_1,\dots, e_n$ satisfying the similar axioms of Theorem  \ref{Theosat}. The models of such a theory are called {\em protomodular algebras}. Protomodular categories  include all
Abelian categories, the category of all groups, loops or even semi-loops, rings with or without unit, associative algebras with or without
unit, Lie algebras, Jordan algebras, Boolean algebras, Heyting algebras, Boolean rings, Heyting semi-lattices, and so on.  If in Definition \ref{TopologicalTalgebra},  $\T$ is a protomodular theory, the corresponding topological $\T$-algebras are  called {\em topological protomodular algebras}.
It was proved in \cite{Bor-Man2} that some results about  topological semi-abelian algebras studied in \cite{Bor-Man} can be extended to the topological protomodular algebras.

We also recall that in an algebraic theory $\T$, the set of constants  is the free algebra on the empty set of generators and is trivially the initial object in the category $\C$ of $\T$-algebras.  If $\T$ has a  unique constant $e$, the initial object is thus reduced to the singleton $\{e\}$ and therefore, becomes trivially also a final object, that is, a zero object in  $\C$. A  theory $\T$ is equivalent to the dual of the category of finitely generated $\T$-algebras. In this equivalence, the object $T^n$ of the theory corresponds by duality to the free algebra $F(n)$ on $n $ generators. In particular the object $T^0$ corresponds to the free algebra on the empty set, that is, the zero algebra $\{e\}$.
We  can now prove the following  Lemma which is used later in some proofs.
\begin{Lem}\label{Lemmapsto-e} Let $\T$ be a semi-abelian theory and $A$ a $\T$-algebra with constant $e\in A$. Then an $n$-ary mapping $\tau\colon A^n\rightarrow A$ maps $(e,\dots,e)$ to $e$. \end{Lem}
\begin{Prf} Let  $\tau\colon A^n \rightarrow A$ be an $n$-ary mapping of the theory $\T$. It corresponds by duality to a $\T$-homomorphism
$t\colon  F(1) \rightarrow  F(n)$ between the corresponding free algebras. By the definition of a free algebra, given a $\T$-algebra $A$, we have a bijection between
\[    \Set(\{1,2,\dots,n\},A) \rightarrow \C(F(n),A) \]
Choosing $n$-elements $a_1,\dots,a_n$ in the $\T$-algebra $A$ is the same as choosing a $\T$-homomorphism $\beta\colon  F(n) \rightarrow A$.  The composite
                $\beta \circ  t \colon  F(1)\rightarrow  F(n) \rightarrow A$ of $t$ and $\beta$ is a $\T$-homomorphism $F(1)\rightarrow  A$ and corresponds thus to the  single element $\tau(a_1,\dots,a_n)$ of $A$.

Suppose now that $a_1$ to $a_n$ are all $e$. This means that the $\T$-homomorphism $\beta$ maps all generators of $F(n)$ to $e$, that is, $\beta$ factors through the zero object $\{e\}$

                \[\beta\colon F(n)\rightarrow \{e\} \rightarrow A\]

But then the composite $T$-homomorphism $\beta \circ t$ factors through $\{e\}$ as well

               \[ \beta \circ  t \colon  F(1) \rightarrow F(n) \rightarrow \{e\}\rightarrow  A\]

This composite is thus the zero $T$-homomorphism, which maps the generator of $F(1)$ on $e$ in $A$. Hence we have that $\tau(e,\dots,e)=e$.

\end{Prf}

We remind the following construction from \cite[pp.295-302]{Rot}.

Let $X$ be a topological space with a base point $x_0$ and  $G$ a subgroup of the fundamental group $\pi_1(X,x_0)$. Let   $\P(X,x_0)$ be  the
set of all paths  $\beta$ in $X$ with source point $x_0$. Then the relation   on $\P(X,x_0)$ defined by $\beta\simeq \gamma$ if and only if $\beta(1)=\gamma(1)$ and $[\beta\circ \gamma^{-1}]\in G$,  is an equivalence relation, where `$\circ$' denotes the concatenation of the paths.  Denote the
equivalence class of $\beta$ by $\langle \beta\rangle _G$ and
define $\widetilde{X}_G$ as the set of all such equivalence
classes of the paths in  $X$ with source point $x_0$. Define a function $p\colon \widetilde{X}_G\rightarrow
X$ by $p(\langle \beta\rangle _G)=\beta(1)$. Let $\beta_0$ be
the constant path at $x_0$ and $\tilde{x}_0=\langle
\beta_0\rangle _G\in \widetilde{X}_G$. If  $\beta\in \P(X,x_0)$ and $U$ is  an open
neighbourhood of $\beta(1)$, then  a path  of  the form
$\beta\circ \lambda$, where $\lambda$ is a path in $U$ with
$\lambda(0)=\beta(1)$, is called a {\em continuation} of
$\beta$. For a $\langle \beta\rangle _G\in \widetilde{X}_G$ and an open
neighbourhood $U$ of $\beta(1)$, let $(\langle \beta\rangle _G,U)=\{\langle \beta\circ \lambda \rangle_G\colon ~\lambda(I)\subseteq
U\}$. Then the subsets
$(\langle \beta\rangle _G, U)$ form a basis for a topology on
$\widetilde{X}_G$ such that the map $p\colon
(\widetilde{X}_G,\tilde{x}_0)\rightarrow (X,x_0)$ is
continuous \cite[Lemma 10.31]{Rot}. We also know from  \cite[Theorem 10.34]{Rot} that if $X$ is connected and has a universal cover, then $p\colon
(\widetilde{X}_G,\tilde{x}_0)\rightarrow (X,x_0)$ is a
covering map with characteristic  group $G$.

\section{Topological $\T$-algebras}
In this section for a semi-abelian theory $\T$,  we apply the construction of  \cite[pp.295-302]{Rot} stated above to the topological $\T$-algebras.  We first need the following preparation:

Let $\T$ be a semi-abelian theory, $A$ a topological $\T$-algebra with the unique constant  $e$ and  $\P(A,e)$  the set of all paths in $A$ with source points $e$. For every $n$-ary continuous mapping $\tau \colon A^n\rightarrow A$ and the paths  $\beta_1,\dots, \beta_{n}$ of $\P(A,e)$   we have a continuous  mapping $[0,1]\rightarrow A$ defined by
\begin{align}\label{thaucomposite}
\tau(\beta_1,\dots,\beta_{n})(t)&=\tau(\beta_1(t),\dots, \beta_{n}(t))~~~
    \end{align}
for $t\in [0,1]$.
Then by Lemma \ref{Lemmapsto-e} we have that
\begin{align*}
\tau(\beta_1,\dots, \beta_{n})(0)&= \tau(\beta_1(0),\dots, \beta_{n}(0))\\
                &=\tau(e,\dots, e) \\
               &=e  \\
\end{align*}
and therefore  $\tau(\beta_1,\dots, \beta_{n})$ is a path of $\P(A,e)$.

We also have  \begin{align}\label{abnininvers}
(\tau(\beta_1,\dots, \beta_{n}))^{-1}= \tau({\beta_1}^{-1},\dots, {\beta_{n}}^{-1})\end{align}
where, for a path , say $\beta$ the inverse $\beta^{-1}$ denotes the inverse path  defined by $\beta^{-1}(t)=\beta(1-t)$ for $t\in [0,1]$.  Then by the evaluation of the concatenation  of the paths in  $A$    at $t\in [0,1]$     we have that  the  {\em interchange rule}
      \begin{align}\label{thetecomposite}
\tau(\beta_1,\dots,\beta_{n})\circ \tau (\gamma_1,\dots, \gamma_{n})=
\tau(\beta_1\circ \gamma_1,\dots, \beta_{n}\circ \gamma_{n})
    \end{align}
holds whenever the concatenations of the paths $\beta_i$ and $\gamma_i$  are defined, where  `$\circ$'  denotes the concatenation of the paths.
More precisely evaluating the concatenations of these paths  at $t\in [0,1]$  for the left side of Eq.\ref{thetecomposite} we have

\[(\tau(\beta_1,\dots,\beta_{n})\circ \tau (\gamma_1,\dots, \gamma_{n}))(t) =\left\{\begin{array}{ll}
                \tau(\beta_1(2t),\dots, \beta_{n}(2t)),  &  \mbox {$0\leq t\leq \frac{1}{2}$}\\\\
          \tau(\gamma_1(2t-1),\dots, \gamma_{n}(2t-1)),   &  \mbox{$ \frac{1}{2}\leq t\leq 1$}
                \end{array}
                \right. \]
and for the right side
\begin{align*}
(\tau(\beta_1\circ \gamma_1,\dots, \beta_{n}\circ \gamma_{n})(t))&= \tau((\beta_1\circ \gamma_1)(t),\dots, (\beta_{n}\circ \gamma_{n}(t))\\
                            &=\tau(\beta_1(2t),\dots, \beta_{n}(2t))
               \end{align*}
if  $ 0\leq t\leq \frac{1}{2}$
and
\begin{align*}
\tau(\beta_1\circ \gamma_1,\dots, \beta_{n}\circ \gamma_{n})(t)&= \tau((\beta_1\circ \gamma_1)(t),\dots, (\beta_{n}\circ \gamma_{n}(t))\\
                            &=\tau(\gamma_1(2t-1),\dots, \gamma_{n}(2t-1))
               \end{align*}
if for $ \frac{1}{2} \leq t\leq 1$.
This proves that Eq.\ref{thetecomposite} is satisfied.

\begin{Thm}\label{FundgroupTalgebra} Let $\T$ be a semi-abelian theory. If  $A$ is a topological  $\T$-algebra with the unique constant  $e$,  then the fundamental group $\pi_1(A,e)$ becomes a $\T$-algebra.\end{Thm}
\begin{Prf} Let $A$ be a topological $\T$-algebra  with the constant element $e$. Hence categorically it represents a  product preserving functor $F_A\colon \T\rightarrow \Top$ in which $F_A(T)=A$,  for the distinguish object $T$ of $\T$.
Then  for an $n$-ary operation $\tau\colon T^n\rightarrow T$ of the theory $\T$ we have  a  mapping
\[ \pi_1(A,e)^n\rightarrow \pi_1(A,e)\] defined by
\begin{align} \label{definitionTheta}
([\beta_1],\dots,[\beta_{n}])\mapsto [\tau(\beta_1,\dots,\beta_{n})]
\end{align}
for   $[\beta_i]\in  \pi_1(A,e)$ ($1\leq i\leq n$). Here $\tau(\beta_1,\dots,\beta_{n})$ is the path defined by Eq.\ref{thaucomposite}.  The mapping defined on $\pi_1(A,e)^n$ is  well defined by  the continuity of the mapping $\tau\colon A^n\rightarrow A$.   We now prove that according to these mappings,  $\pi_1(A,e)$ becomes a  $\T$-algebra with the constant element $\tilde{e}$, which is the homotopy class of the constant path at  $e\in A$.
An arrow \[(\tau_0,\dots,\tau_{m-1})\colon T^n\rightarrow T^m\] of the theory $\T$  constitutes the mapping \[ \pi_1(A,e)^n\rightarrow \pi_1(A,e)^m\] defined by
  \[
   ([\beta_1],\dots,[\beta_{n}])\mapsto( [\tau_0(\beta_1,\dots,\beta_{n})],\dots,[\tau_{m-1}(\beta_1,\dots, \beta_n)]
\]
Hence by this evaluation we have  a product preserving functor
\[\pi_1 F_A\colon \T\rightarrow \Set\] induced by the distinguished object $T$ of the theory $\T$  \[(\pi_1 F_A)(T)=\pi_1((F_A(T),e)=\pi_1(A,e)\]
The axioms of Theorem \ref{Theosat} can be checked as follows.
\[\alpha_i([\beta],[\beta])=[\alpha_i(\beta,\beta)]=\tilde{e}\] for any  binary mapping $\alpha_{i}$ and $n+1$-ary mapping $\theta$
\begin{align*}
\theta(\alpha_1([\beta],[\gamma]), \dots, \alpha_n([\beta],[\gamma]), [\gamma])&=\theta([\alpha_1(\beta,\gamma)],\dots [\alpha_n(\beta,\gamma)],[\gamma])  \tag{by  Eq.\ref{definitionTheta}}\\
&=[\theta(\alpha_1(\beta,\gamma),\dots \alpha_n(\beta,\gamma),\gamma)] \tag{by  Eq.\ref{definitionTheta}} \\
                            &=[\beta] \tag{by Theorem \ref{Theosat}}
               \end{align*}
for $[\beta],[\gamma]\in \pi_1(A,e)$.  Therefore $\pi_1(A,e)$ becomes a $\T$-algebra for the same semi-abelian theory $\T$.
\end{Prf}

\begin{Prop}\label{Propfunctor1} Let $\T$ be a semi-abelian theory. Then we have  a functor $\pi_1\colon \Top^{\T}\rightarrow \C$ assigning each topological $\T$-algebra $A$ to the $\T$-algebra $\pi_1(A,e)$.
\end{Prop}
\begin{Prf} Let $f\colon A\rightarrow B$ be  a continuous $\T$-homomorphism.  Then by the following evaluation,  the induced map $f_{\star}=\pi_1 f\colon \pi_1(A,e)\rightarrow \pi_1(B,e)$ becomes a $\T$-homomorphism
\begin{align*}
f_{\star}(\tau ([\gamma_1],\dots, [\gamma_n])&=f_{\star}[\tau (\gamma_1,\dots, \gamma_n)]\\
                                            &=[f(\tau (\gamma_1,\dots, \gamma_n))]\\
                                             &=[\tau (f\gamma_1,\dots, f\gamma_n)]\\
                                              &=\tau([ f\gamma_1],\dots, [f\gamma_n])\\
                                              &=\tau(f_{\star}[\gamma_1],\dots, f_{\star}[\gamma_n])
\end{align*}
The axioms for $\pi_1$ to be a functor are straightforward and therefore omitted.
\end{Prf}

We recall from \cite[Theorem A.2]{Bor-Man} that for a semi-abelian theory $\T$,  a sub $\T$-algebra $B$ of a  $\T$-algebra $A$ is {\em normal} if for every operation $\tau (X_1,\dots, X_k, Y_1,\dots, Y_l)$ of the theory such that $\tau(X_1,\dots,X_k,e,\dots, e)=e$ one has $\tau(a_1,\dots,a_k,b_1,\dots, b_l)\in B$ for $a_1,\dots,a_k\in A$ and $b_1,\dots,b_l\in B$.

\begin{Lem} \label{Fundgpnormal} Let $A$ be a topological $\T$-algebra for a semi-abelian theory $\T$ with constant $e$. If $B$ is a normal sub $\T$-algebra of $A$, then $\pi_1(B,e)$ becomes a normal sub $\T$-algebra of $\pi_1(A,e)$.
\end{Lem}
\begin{Prf} Let $\tau (X_1,\dots, X_k, Y_1,\dots, Y_l)$  be an operation of the theory such that $\tau(X_1,\dots,X_k,e,\dots, e)=e$; and let  $[\gamma_1],\dots,[\gamma_k]\in \pi_1(A,e)$ and $[\beta_1],\dots,[\beta_n]\in \pi_1(B,e)$. Since $B$ is a normal sub $\T$-algebra of $A$ we have
\[\tau(\gamma_1,\dots, \gamma_k, \beta_1, \dots, \beta_l)(t)=\tau(\gamma_1(t),\dots, \gamma_k(t), \beta_1(t), \dots, \beta_l(t))\in B\]
and therefore
\begin{align*}
\tau ([\gamma_1],\dots, [\gamma_k], [\beta_1], \dots, [\beta_l])&=[\tau(\gamma_1,\dots, \gamma_k, \beta_1, \dots, \beta_l)] \in \pi_1(B,e)
\end{align*}
Hence $\pi_1(B,e)$ is a normal sub $\T$-algebra of $\pi_1(A,e)$.
\end{Prf}

Let $\mathbb{T}$ be a semi-abelian theory with the constant $e$ and $A$  a  topological $\T$-algebra. Assume that $A'$ is a sub $\T$-algebra of $A$ and $B$ is the fundamental group $\pi_1(A',e)$. Then  $B$ is a sub $\T$-algebra of $\pi_1(A,e)$; and therefore by  \cite[Lemma 10.31]{Rot}  we have a topological space  $\widetilde{A}_B$ and a continuous map $p\colon (\widetilde{A}_B,\tilde{e})\rightarrow (A,e)$  between topological spaces, corresponding to the subgroup $B$ of $\pi_1(A,e)$. Hence  we can prove the following theorem for topological $\T$-algebras.

\begin{Thm}\label{CovCorrG}  Let $\mathbb{T}$ be a semi-abelian theory with the constant $e$ and $A$  a  topological $\T$-algebra. Let  $A'$ be a sub $\T$-algebra of $A$ and $B$  the fundamental group $\pi_1(A',e)$.  Then  $\widetilde{A}_B$ becomes a topological  $\T$-algebra such that  the map $p\colon (\widetilde{A}_B,\tilde{e})\rightarrow (A,e)$  is a continuous $\T$-homomorphism.\end{Thm}
\begin{Prf} Since $B$ is a subgroup of the fundamental group $\pi_1(A,e)$, by \cite[Lemma 10.31]{Rot} we have a  continuous map  $p\colon (\widetilde{A}_B,\tilde{e})\rightarrow (A,e)$ of  topological spaces corresponding to  $B$.  Hence    $\widetilde{A}_B$ is defined as the set of equivalence classes defined via $B$.   Then each    mapping   defined by Eq.\ref{thaucomposite}  reduces a mapping  defined by
   \begin{align}\label{Operationtheta}
    \tau(\langle \beta_1\rangle,\cdots,\langle \beta_{n}\rangle)= \langle\tau(\beta_1,\dots,\beta_{n})\rangle.
\end{align}
for $\langle \beta_1\rangle,\cdots,\langle \beta_{n}\rangle\in \widetilde{A}_B$.
We now prove that this  mapping is well defined:
For the paths $\beta_1,\dots\dots,\beta_{n}$ and $\gamma_1,\dots,\gamma_{n}$ in $\P(A,e)$ with $\beta_i(1)=\gamma_i(1)$ one has that
  \begin{align*}
  [\tau(\beta_1,\dots,\beta_{n}) \circ(\tau(\gamma_1,\dots,\gamma_{n}))^{-1}]&= [\tau(\beta_1,\dots,\beta_{n})\circ \tau({\gamma_1}^{-1},\dots,{\gamma_{n}}^{-1})] \tag{by  Eq.\ref{abnininvers}}\\
  &=[( \tau(\beta_1 \circ{\gamma_1}^{-1},\dots,\beta_{n}\circ {\gamma_{n}}^{-1})] \tag{by  Eq.\ref{thetecomposite}}\\
   &=\tau([\beta_1 \circ{\gamma_1}^{-1}],\dots,[\beta_{n}\circ {\gamma_{n}}^{-1}]) \tag{by Eq.\ref{definitionTheta}}
\end{align*}
Since $B$ is a sub $\T$-algebra of $\pi_1(A,e)$ one has   $\tau([\beta_1 \circ{\gamma_1}^{-1}],\dots,[\beta_{n} \circ{\gamma_{n}}^{-1}])\in B$ when  $[\beta_i\circ\gamma_i^{-1}]\in B$ $(1\leq i\leq n)$.  Hence $[\tau(\beta_1,\dots,\beta_{n}) \circ(\tau(\gamma_1,\dots,\gamma_{n}))^{-1})]\in B$ whenever $[\beta_i\circ \gamma_i^{-1}]\in B$  and therefore the mapping $\tau$ defined on $\widetilde{A}_B$ is well defined.  Hence we have a product preserving functor  $F_B\colon \T\rightarrow \Set$ induced  by $F_B(T)=\widetilde{A}_B$, for the distinguished object of $\T$. The axioms of Theorem \ref{Theosat} for the  mappings defined in Eq.\ref{Operationtheta} are satisfied and hence  $\widetilde{A}_B$ becomes a semi-abelian algebra. The map  $p$ is a  $\T$-homomorphism by the details
\begin{align*}
  p(\tau(\langle \beta_1\rangle,\cdots,\langle \beta_{n}\rangle))&= p(\langle\tau(\beta_1,\dots,\beta_{n})\rangle)\tag{by  Eq.\ref{Operationtheta}}\\
  & =(\tau(\beta_1,\dots,\beta_{n}))(1) \\
  &=\tau(\beta_1(1),\dots,\beta_{n}(1))  \tag{by  Eq.\ref{thaucomposite}}  \\
  & =\tau(p(\langle \beta_1\rangle),\cdots,p(\langle \beta_{n}\rangle))\end{align*}
  for $\langle \beta_1\rangle,\cdots,\langle \beta_{n}\rangle\in \widetilde{A}_B$.

To prove that $\widetilde{A}_B$ is a topological $\T$-algebra, we now prove that each  $n$-ary mapping $\tau$ of $\widetilde{A}_B$ defined by Eq.\ref{Operationtheta} is  continuous.  Let $\beta=(\beta_1,\dots,\beta_{n})$  and let $(V, \langle\tau(\beta)\rangle )$ be a base open  neighbourhood of  $\langle\tau(\beta)\rangle $. Then   $V$ is an open neighbourhood of \[\tau(\beta(1))=\tau(\beta_1(1),\dots,\beta_{n}(1))\] and since the mapping $\tau$ on $A$ is continuous, there are respectively open neighbourhoods $U_1,\dots,U_{n}$ of $\beta_1(1),\dots,\beta_n(1)$  such that \[\tau(U_1\times\dots \times U_{n})\subseteq V.\] Setting $U=U_1\times\dots \times U_{n}$, one obtains
$\tau(U,\langle \beta\rangle )\subseteq (V,\langle\tau (\beta)\rangle )$  which concludes that $n$-ary mapping $\tau'$ on  $\widetilde{A}_B$ is continuous (see the proof of \cite[Lemma 10.31]{Rot}).
  \end{Prf}

  \section{Covering morphisms of internal groupoids in semi-abelian categories}

  Let $\T$ be a semi-abelian theory.  In this section we define internal groupoid  in the semi-abelian category  $\C$ of $\T$-algebras and extend some results of  \cite{Ak-Na-Mu-Tu} about the coverings of  internal groupoids in the groups with operations to  the  internal groupoids in the category $\C$.  To define an internal groupoid  in the semi-abelian category  $\C$ of $\T$-algebras we comply with the notations for groupoids given in  Section \ref{PrelimCov}.

Similar to the notion of an internal category in the category of groups with operations as defined in \cite{Por}, the definition of internal groupoid in the category $\C$ of models for a semi-abelian theory $\T$ is given as follows.
\begin{Def}\label{Internalgpd} {\em Let $\T$ be a semi-abelian theory.  An {\em internal groupoid} in $\C$  is a groupoid  $G$ in which  the set $\Ob(G)$ of objects  and the set $G$ of arrows  are both  $\T$-algebras; and the source and target point maps $s,t\colon G\rightarrow \Ob(G)$, the object inclusion map  $\epsilon\colon \Ob(G)\rightarrow G$, the partial composite $\circ\colon G_{t}\times_{s}
G\rightarrow G,(g,h)\mapsto g\circ h$ and the inversion $G\rightarrow G,g\mapsto g^{-1}$ are all  $\T$-homomorphisms.\qed }\end{Def}

Note  that the partial composite `$\circ$' is a $\T$-homomorphism if and only if for every $n$-ary mapping $\tau$  the  {\em interchange rule}
\begin{align}\label{intechangeint}
\tau(g_1\circ h_1,\dots,g_{n}\circ h_{n})=\tau(g_1,\dots,g_{n})\circ \tau(h_1,\dots,h_{n})
    \end{align}
is satisfied  for $g_1,\dots,g_n \in G$ and  $h_1,\dots,h_n\in G$ whenever  one side composition is defined.
For the category of internal groupoids  in $\C$ we use  the  notation $\Cat(\C)$.

In particular if $\T$ is the group theory, then  $\C$ becomes  the  category of groups and hence  an internal groupoid in $\C$ becomes  a  {\em group-groupoid} which is also called in literature  as {\em group objects} or 2-{\em group}.

\begin{Rem}\label{inverseofmappings}{\em  Let $G$ be an internal groupoid in $\C$  for a semi-abelian theory $\T$. Then we have the following:

\begin{enumerate}
\item Since the inversion $G\rightarrow G,g\mapsto g^{-1}$ is a $\T$-homomorphism for any $n$-ary mapping $\tau$  and $g_1,\dots,g_n\in G$ we have
  \begin{align}\label{invermapping}
(\tau(g_1,\dots, g_n))^{-1}=\tau({g_1}^{-1},\dots,{g_n}^{-1})
    \end{align}

\item  Since the object inclusion map  $\epsilon\colon \Ob(G)\rightarrow G, x\rightarrow  1_x$ is a $\T$-homomorphism, the identity arrow $1_e$ at the constant $e\in \Ob(G)$ is the constant arrow of $G$ and
    \begin{align*}
1_{\tau(x_1,\dots, x_n)=\tau(1_{x_1},\dots,1_{x_n})}
    \end{align*}

\end{enumerate}
}\end{Rem}

\begin{Lem}\label{StGT-algebra} For a semi-abelian theory $\T$, if  $G$ is  an internal groupoid in $\C$ with unique constant $e\in \Ob(G)$, then we have the following:
\begin{enumerate}

\item $\St_Ge$ is a sub $\T$-algebra of $G$.

\item   $G(e)$, the object group at $e\in \Ob(G)$, is a sub $\T$-algebra of $G$.
\end{enumerate}
\end{Lem}
\begin{Prf}
\begin{enumerate}
\item Let $G$ be an internal groupoid in $\C$. Since the source point map $s\colon G\rightarrow \Ob(G)$ is a $\T$-homomorphism, by Lemma \ref{Lemmapsto-e} for $g_1,\dots, g_{n}\in \St_Ge$  we have
\begin{align*}
  s(\tau(g_1,\dots,g_{n}))&=\tau(s(g_1),\dots,s(g_{n}))\\
                                      &=\tau(e,\dots,e)\\
                                      &=e
  \end{align*}
and therefore $\tau(g_1,\dots,g_{n})\in \St_Ge$ when $g_1,\dots,g_{n}\in \St_Ge$. Hence $\St_Ge$ becomes a sub $\T$-algebra of $G$ with unique constant $1_e$.

\item  In addition to the proof of (1) we need to prove that similar axiom are satisfied for the target point map $t\colon G\rightarrow \Ob(G)$.
Since $t$ is a $\T$-homomorphism by Lemma \ref{Lemmapsto-e} for  $g_1,\dots, g_{n}\in G(e)$   we have
\begin{align*}
  t(\tau(g_1,\dots,g_{n}))&=\tau(t(g_1),\dots,t(g_{n}))\\
                                      &=\tau(e,\dots,e)\\
                                      &=e
  \end{align*}
and therefore  $\tau(g_1,\dots,g_{n})\in G(e)$ when $g_1,\dots,g_n\in G(e)$. Hence $G(e)$ becomes a sub $\T$-algebra of $G$.

\end{enumerate}

\end{Prf}

\begin{Exam}{\em Let $\T$ be a semi-abelian theory and let  $A$ be a $\T$-algebra.   Then the groupoid $G=A\times A$ with object set  $A$ such that a pair $(a,b)$ is an arrow from $a$ to $b$  with inverse arrow $(b,a)$ and the  composition is defined by $(a,b)\circ (b,c)=(a,c)$,  becomes an internal groupoid  in $\C$.

Here  an $n$-ary mapping $\tau$ on $G$ is defined by
\[\tau((a_1,b_1),\dots,(a_n,b_n)))=(\tau(a_1,\dots,a_n),\tau(b_1,\dots,b_n)) \]

One can check that for  $g_i=(a_i,b_i)$, $h_i=(b_i,c_i)$ ($1\leq i\leq n$) the following interchange  rule holds
\[
\tau(g_1\circ h_1,\dots,g_{n}\circ h_{n})=\tau(g_1,\dots,g_{n})\circ \tau(h_1,\dots,h_{n})
\]}
\end{Exam}

The following result enables us to produce more examples of the internal groupoids in $\T$-algebras.

\begin{Thm}\label{Fundgpdint}  Let $\T$ be a semi-abelian theory.  If   $A$ is a topological $\T$-algebra, then the fundamental groupoid $\pi A$ becomes an internal groupoid in the semi-abelian category  $\C$ of $\T$-algebras.
\end{Thm}
\begin{Prf} Let   $A$ be a topological $\T$-algebra for a semi-abelian theory $\T$. As similar to Eq.\ref{definitionTheta} one can  define $n$-ary mappings  for  the fundamental groupoid $\pi A$. Since $A$ can be represented by a product preserving functor $F_A\colon \T\rightarrow \Top$ such that $F_A(T)=A$, for the distinguish object $T$ of the theory $\T$  we have a functor $\pi F_A\colon \T\rightarrow \Set$ induced by $(\pi F_A)(T)=\pi(F_A(T))=\pi A$. If $\tau\colon T^n\rightarrow T^m$ is an operation of the theory $\T$, then we have
$(\pi F_A)(\tau)\colon (\pi(A))^n\rightarrow (\pi(A))^m$, a mapping. Hence   $\pi F_A\colon \T\rightarrow \Set$ becomes a product preserving functor.  We can prove that the axioms of Theorem \ref{Theosat} are satisfied as similar to the proof of Theorem \ref{FundgroupTalgebra}.  Therefore $\pi A$, as the set of arrows is a $\T$-algebra.
The interchange rule
\begin{align}
\tau([\beta_1]\circ [\gamma_1],\dots,[\beta_{n}]\circ [\gamma_{n}])=\tau([\beta_1],\dots,[\beta_{n}])\circ \tau([\gamma_1],\dots,[\gamma_{n}])
    \end{align}
for $\pi A$ can be checked by the following evaluating the concatenation  of  the paths:
\begin{align*}
  \tau([\beta_1]\circ [\gamma_1],\dots, [\beta_n]\circ [\gamma_n])&=\tau([\beta_1 \circ\gamma_1],\dots, [\beta_n\circ \gamma_n])\\
  &=[ \tau(\beta_1 \circ{\gamma_1},\dots,\beta_{n}\circ {\gamma_{n}})] \tag{by Eq.\ref{definitionTheta}}\\
   &=[\tau(\beta_1,\dots, \beta_n)\circ \tau(\gamma_1,\dots, \gamma_n) ]\tag{by  Eq.\ref{thetecomposite}}\\
     &=[\tau(\beta_1,\dots, \beta_n)]\circ [ \tau(\gamma_1,\dots, \gamma_n) ]\\
   &=\tau([\beta_1],\dots, [\beta_n])\circ \tau([\gamma_1],\dots, [\beta_n])\tag{by Eq.\ref{definitionTheta}}
\end{align*}
The other details to complete the proof are straightforward.\end{Prf}

As a result of Theorem \ref{Fundgpdint} for a semi-abelian theory $\T$, we have a   functor $\pi\colon  \Top^{\T} \rightarrow \Cat(\C)$  assigning each topological $\T$-algebra $A$ to the internal groupoid $\pi(A)$ in $\C$.

\begin{Lem} Let $A$ be a topological $\T$-algebra for a semi-abelian theory $\T$. If $B$ is a normal sub $\T$-algebra of $A$, then $\pi(B)$ becomes a normal sub $\T$-algebra of $\pi(A)$.
\end{Lem}
\begin{Prf} The proof can be done as similar to the proof of Lemma \ref{Fundgpnormal}.
\end{Prf}

\begin{Prop} For a   semi-abelian theory $\T$; if $A$ and  $B$ are topological $\T$-algebras, then  $\pi (A\times B)$ and $\pi A \times \pi B$ are  isomorphic as internal groupoids in $\C$.
\end{Prop}
\begin{Prf} By \cite[6.4.4]{Br} we know that the   map $f\colon \pi(A\times B)\rightarrow \pi A\times \pi  B$ defined by  $f([\beta])=([\beta_A],[\beta_B])$ for a homotopy class $[\beta]\in \pi(A\times B)$ is an isomorphism of the underlying groupoids, where $\beta_A$ and $\beta_B$ are the projections of the path $\beta$ on  $A$ and $B$ respectively.

Replacing $A$ with $A\times B$ in Theorem \ref{Fundgpdint} we have that $\pi(A\times B)$ is an internal groupoid, where $n$-ary mappings of $\pi(A\times B)$ are defined  by
\begin{align*}
   {\tau}([\beta_1],\dots,[\beta_{n}])&= [\tau(\beta_1,\dots,\beta_{n})]\\
  &=[\tau(({\beta_1}_A,{\beta_1}_B), \dots,({\beta_n}_A,{\beta_n}_B))]\\
  &=[(\tau({\beta_1}_A,\dots, {\beta_n}_A), \tau({\beta_1}_B,\dots, {\beta_n}_B))]
   \end{align*}
for $[\beta_1],\dots, [\beta_n]\in \pi(A\times B)$.
Further we now  check  that $f$ is a morphism of the internal groupoids in $\C$.
\[   f(\tau([\beta_1],\dots,[\beta_{n}]))=([\tau({\beta_1}_A,\dots, {\beta_n}_A)], [\tau({\beta_1}_B,\dots, {\beta_n}_B)])\]
On the other hand
  \begin{align*}
   {\tau}(f[\beta_1],\dots,f[\beta_{n}])&= \tau(([{\beta_1}_A],[{\beta_1}_B]),\dots, ([{\beta_n}_A],[{\beta_n}_B])) \\
  &=(\tau([{\beta_1}_A],\dots,[{\beta_n}_A]), \tau([{\beta_1}_B],\dots,[{\beta_n}_B]) )\\
   &=[\tau({\beta_1}_A,\dots,{\beta_n}_A), \tau({\beta_1}_B,\dots,{\beta_n}_B)] \\
   \end{align*}
Hence by comparing  these we have that
\[   f(\tau([\beta_1],\dots,[\beta_{n}]))= {\tau}(f[\beta_1],\dots,f[\beta_{n}])\]
and therefore $f$ is a morphism of the internal groupoids in $\C$.
\end{Prf}

Let $G$ be  an internal groupoid  in a certain category  of groups with operations and $X$ a group with operations.
The action of $G$ on $X$ is defined  in \cite[Definition 3.11]{Ak-Na-Mu-Tu}. This definition  is  generalized to the internal groupoids in the category $\C$ of semi-abelian $\T$-algebras as follows:

\begin{Def}{\em  Let $\T$ be a semi-abelian theory and  $G$  an internal groupoid in the category  $\C$ of semi-abelian algebras. Let  $A$ be a $\T$-algebra  and  $\omega\colon A\rightarrow \Ob(G)$ a $\T$-homomorphism.  If the underlying groupoid  $G$ acts on the underlying set of  $A$  via $\omega$ such that $\varphi\colon A_\omega \times _{s} G\rightarrow A, (a,g)\mapsto ag$ is also a $\T$-homomorphism, then we say that the internal groupoid $G$ {\em acts} on $\T$-algebra $A$ via $\omega$.\qed}\end{Def}
We write $(A,\omega,\varphi)$ for such an action. Here note that  $\varphi\colon A_\omega\times _{s} G\rightarrow A, (a,g)\mapsto ag$ is a $\T$-homomorphism in $\C$ if and only if for $a_1,\dots, a_n\in A$ and $g_1,\dots,g_{n}\in G$
  \begin{align}\label{2}
  \tau(a_1,\dots, a_n) \tau(g_1,\dots,g_{n})= \tau(a_1g_1,\dots,a_n g_{n})
    \end{align}
   whenever  one side is defined.

\begin{Exam}\label{Examaction}{\em Let $\T$ be a semi-abelian theory and  $p\colon H\rightarrow G$   a  morphism of internal groupoids in the semi-abelian category $\C$ of $\T$-algebras such that $p$ is a covering morphism on the underlying groupoids. Then the internal groupoid  $G$ acts on the $\T$-algebra $A=\Ob(H)$ via $\Ob(p)\colon A\rightarrow \Ob(G)$ assigning to $a\in A$ and $g\in \St_Gp(a)$ the target of the unique lifting $\widetilde{g}$   of $g$ in $H$ with source $a$. Clearly the underlying groupoid of $G$ acts on the underlying set   and  by evaluating the  uniqueness of the lifting,  the condition Eq.\ref{2}  is satisfied for  $a_1,\dots,a_n\in A$ and  $g_1,\dots,g_n\in G$ whenever one side is defined.\qed}\end{Exam}

The Characterization of Theorem \ref{Theoactiongpdcover} for semi-abelian theory is as follows:

\begin{Thm} \label{semiabelactiongpd} Let $\T$ be a  semi-abelian theory with unique constant $e$ and    $G$ an internal groupoid in $\C$ such that the underlying groupoid is transitive.  Let  $G(e)$  be the object group at $e\in \Ob(G)$  and  $C$ a subgroup and a sub $\T$-algebra   of  $G(e)$. Then the set $A_C$  of  cosets $C\circ g=\{c\circ g\mid c\in\ C\}$ for $g$ in $\St_G e$ becomes a $\T$-algebra and  the internal  groupoid $G$ acts on  $A_C$ by $(C\circ a)g=C\circ a\circ g$.
\end{Thm}
\begin{Prf}  Define $n$-ary mappings on $A_C$,  the set of the cosets $C\circ g$,  by
  \begin{align}\label{definitionThetaonCoa}
    \tau(C\circ g_1,\cdots,C\circ g_{n})= C\circ \tau(g_1,\dots,g_{n})
\end{align}
Here note that by Lemma \ref{StGT-algebra} (1),  $\tau(g_1,\dots,g_{n})\in \St_G e$ whenever $g_1,\dots,g_{n}\in \St_Ge$. We now prove that the $n$-ary mappings  $\tau$ are well defined. Let  $C\circ g_1=C\circ h_1$, $\dots$, and $C\circ g_n=C\circ h_n$. Since $C$ is a subgroup of $G(e)$ we have $g_1 \circ {h_1}^{-1}\in C$,   $\dots$, $g_{n} \circ {h_{n}}^{-1}\in C$ and therefore
    \begin{align*}
  \tau(g_1,\dots,g_{n}) \circ (\tau(h_1,\dots,h_{n}))^{-1}&= \tau(g_1,\dots,g_{n})\circ \tau({h_1}^{-1},\dots,{h_{n}}^{-1}) \tag{by Eq.\ref{invermapping}}\\
  &= \tau(g_1 \circ {h_1}^{-1},\dots,g_{n}\circ {h_{n}}^{-1}) \tag{by  Eq.\ref{intechangeint}}\\
  \end{align*}
Since $C$ is a sub   $\T$-algebra of $G(e)$ we have that $\tau (g_1 \circ {h_1}^{-1},\dots,g_{n} \circ {h_{n}}^{-1})\in C$ and so  $\tau(g_1,\dots,g_{n}) \circ (\tau(h_1,\dots,h_{n}))^{-1}\in C$. Hence $C\circ \tau(g_1,\dots,g_{n})=C\circ\tau(h_1,\dots,h_{n})$ and the $n$-ary mappings $\tau$ defined by Eq.\ref{definitionThetaonCoa} are well defined. Hence we have a functor $F_C\colon \T\rightarrow \Set$ defined by $F_C(T)=A_C$.
Moreover for $2$-ary mappings  $\alpha_i$ by Eq.\ref{definitionThetaonCoa} we have \[\alpha_i(C\circ g,C\circ g)=C\circ \alpha_i(g,g)=C\circ e=C\]
By the Eq.\ref{definitionThetaonCoa} we have the   following evaluation.
   \begin{align*}
  \theta(\alpha_1(C\circ g,C\circ h),\dots, \alpha_n(C\circ g,C\circ h),C\circ h)&=\theta(C\circ \alpha_1(g,h), \dots, C\circ \alpha_n(g,h),C\circ h) \\
    &=C\circ \theta(\alpha_1(g,h),\dots, \alpha_n(g,h), h)\\
   &=C\circ g \tag{by Theorem \ref{Theosat} (\ref{sattheta}) }
   \end{align*}
Hence  the axioms of Theorem \ref{Theosat} are satisfied.
Therefore $A_C$ becomes a semi-abelian algebra. Here the underlying groupoid of $G$ acts on the set $A=A_C$ by
  \begin{align}\label{2xx}
   A_{\omega} \times_s G\rightarrow A, (C\circ g,h)\mapsto (C\circ g)h=C\circ g\circ h
   \end{align}
via the map $\omega\colon A_C\rightarrow \Ob(G), C\circ g\mapsto t(g)$.

For $a_i=C\circ g_i\in A$ $(i=1,\dots, n)$ the following evaluations prove that the interchange rule Eq.\ref{2} is satisfied:
\begin{align*}
  \tau(a_1,\dots,a_n) \tau(h_1,,\dots, h_n)&=\tau(C\circ g_1,\dots,C\circ g_n) \tau(h_1,\dots, h_n)\\
  &=(C\circ \tau(g_1,\dots, g_n))\tau(h_1,\dots, h_n)\tag{by Eq.\ref{definitionThetaonCoa}}\\
  &=C\circ \tau(g_1, \dots, g_n)\circ \tau(h_1,\dots, h_n)\tag{by Eq.\ref{2xx}}\\
   &=C\circ \tau(g_1\circ h_1,\dots, g_n\circ h_n)\tag{by Eq.\ref{intechangeint}}\\
   \end{align*}
and
\begin{align*}
  \tau(a_1h_1,\cdot \dots, a_nh_n)&= \tau((C\circ g_1)h_1,\dots, (C\circ g_n)h_n)\\
  &=\tau(C\circ g_1\circ h_1,\dots, (C\circ g_n \circ h_n)\tag{by Eq.\ref{2xx}}\\
   &=C\circ \tau(g_1\circ h_1,\dots,  g_n\circ h_n)\tag{by Eq.\ref{definitionThetaonCoa}}
     \end{align*}
Hence $\tau(a_1,\dots,a_n) \tau(h_1,,\dots, h_n)= \tau(a_1h_1,\cdot \dots, a_nh_n)$ and therefore the interchange rule in  Eq.\ref{2} is satisfied. \end{Prf}

We now give a concrete  example to Theorem \ref{semiabelactiongpd}.
\begin{Exam} \rm   Let $\T$ be a semi-abelian algebra with a constant $e$ and $A$  a topological $\T$-algebra. Then by Theorem \ref{Fundgpdint},   $G=\pi A$ becomes  an internal groupoid in $\C$ and for a  sub $\T$-algebra $B$ of $A$, the fundamental group $C=\pi_1(B,e)$ becomes a sub $\T$-algebra of $G(e)=\pi_1(A,e)$.  Hence by Theorem, \ref{semiabelactiongpd} the set  $A_C$ of all cosets $\{C\circ g\mid g\in \St_{\pi A}e \} $ is a $\T$-algebra and the internal groupoid  $\pi A$ acts on $A_C$.

\end{Exam}
\begin{Thm}\label{Teointernalgpd} Let  $\T$ be a semi-abelian theory.  Let $G$ be an internal groupoid in $\C$ and $A$ a $\T$-algebra. Suppose that $G$ acts on the $\T$-algebra $A$ via a $\T$-homomorphism $\omega\colon A\rightarrow \Ob(G)$. Then the semi-direct product groupoid $G\ltimes A$   becomes  an internal groupoid in $\C$ such that the projection $p\colon G\ltimes A\rightarrow G$  defined on objects by $\omega$   and on arrows by $(g,a)\mapsto g$  is a morphism of internal groupoids which is a covering morphism on the underlying groupoids.\end{Thm}
\begin{Prf} By \cite[10.42]{Br} we know that the projection map $p\colon G\ltimes A\rightarrow G$ is a covering morphism of groupoids. Then the   semi-direct product  groupoid $H=G\ltimes A$ becomes a $\T$-algebra  by the $n$-ary mappings defined by
\begin{align}\label{defnofalphaiforsemidrecttheta}
\tau((g_1,a_1),\dots(g_{n},a_{n}))=(\tau(g_1,\dots,g_{n}),\tau(a_1,\dots,a_{n}))
\end{align}
In addition the source and target point maps $s,t\colon H\rightarrow
A$, the object inclusion map  $\epsilon\colon A\rightarrow H$
and the partial composition $\circ\colon H_{t}\times_{s}
H\rightarrow H,(h,k)\mapsto h\circ k$ are $\T$-homomorphisms.
Hence $G\ltimes A$ becomes an internal groupoid in $\C$. Moreover for given arrows $(g_1,a_1),\dots, (g_n,a_n)$ of $G\ltimes A$ and any $n$-ary mapping $\tau$   by the following evaluation, $p$ is a $\T$-homomorphism.
\begin{align*}
  p(\tau((g_1,a_1),\dots,(g_n,a_n)))&= p((\tau(g_1,\dots,g_n),\tau(a_1,\dots,a_n)))\\
  &=\tau(g_1,\dots,g_n)\\
   &= \tau(p(g_1,a_1),\dots,p(g_n,a_n))
     \end{align*}
\end{Prf}

Let $G$ be an internal groupoid in $\C$ for a semi-abelian theory $\T$.  Then we have a category  $\Act_{\Cat(\C)}/G$  whose objects are  actions  $(A,\omega,\varphi)$ of the internal groupoid $G$  on $\T$-algebras and  morphisms, say from  $(A,\omega,\varphi)$ to  $(A',\omega',\varphi')$ are $\T$-homomorphisms   $f\colon A\rightarrow  A'$ such that $\omega=\omega'f$ and $f(ag)=(fa)g$ whenever $ag$ is defined;

Let ${\Cov}_{\Cat(\C)}/G$ be the category  whose objects are the morphisms  $p\colon H\rightarrow G$  of internal groupoids in  $\C$ such that $p$  is a covering morphism of groupoids and arrows are commutative diagrams of morphisms of internal groupoids
\begin{center}
\begin{tikzpicture}
\matrix(a)[matrix of math nodes, row sep=3em, column sep=3.5em, text height=1.5ex, text depth=0.25ex]
{H & K  \\  & G \\};
\path[->,font=\scriptsize](a-1-1) edge node[above]{$f$} (a-1-2);
\path[->,font=\scriptsize](a-1-1) edge node[below,left]{$p$} (a-2-2);
\path[->,font=\scriptsize](a-1-2) edge node[right]{$q$} (a-2-2);
\end{tikzpicture}
\end{center}
where $p$ and $q$ are covering morphisms  on the underlying groupoids. Similarly writing such a diagram as a triple $(f;p,q)$, the composition of the arrows in  ${\Cov}_{\Cat(\C)}/G$ is defined by $(f;p,q)\circ (g;q,r)=(fg;p,r)$.

  We can now prove the equivalence of these categories as follows.

\begin{Thm} \label{Covandactareequval} Let $\T$ be a semi-abelian theory and $G$ an  internal groupoid in $\C$. Then the categories $\Act_{\Cat(\C)}/G$ and
${\Cov}_{\Cat(\C)}/G$  are equivalent.\end{Thm}
\begin{Prf} If $(A,\omega,\varphi)$ is an object of $\Act_{\Cat(\C)}/G$, then by Theorem \ref{Teointernalgpd}, we have a
morphism $p\colon G\ltimes A\rightarrow G$ of internal groupoids in $\C$, which is a covering morphism on the underlying groupoids.  This gives us a functor \[\Gamma\colon
\Act_{\Cat(\C)}/G\rightarrow {\Cov}_{\Cat(\C)}/G.\]

Conversely if  $p\colon H\rightarrow G$ is a  morphism of internal groupoids in $\C$ which is a covering morphism on the underlying groupoids, then by Example \ref{Examaction} we have an action of the internal groupoid $G$ on the $\T$-algebra $A=\Ob(H)$ via $p\colon A\rightarrow \Ob(H)$.
In this way we define a functor  \[\Phi\colon
{\Cov}_{\Cat(\C)}/G\rightarrow \Act_{\Cat(\C)}/G.\]

The natural equivalences $\Gamma\Phi\simeq 1$ and $\Phi\Gamma\simeq 1$ follow.
\end{Prf}

\begin{Def}{\em Let $\T$ be a semi-abelian theory, $H$  a groupoid and   $G$ an internal  groupoid in $\C$ with constant  $e\in \Ob(G)$.  Suppose that   $p\colon H\rightarrow G$ is a covering morphism of groupoids and $e'\in \Ob(H)$.  We say that $\T$-algebraic structure of  $G$ {\em lifts} to $H$ if $H$ becomes an internal groupoid in $\C$ such that $p$ is a morphism of internal groupoids.\qed}\end{Def}

Using Theorem  \ref{semiabelactiongpd} and  Theorem \ref{Teointernalgpd}  we now give a criterion for  the $\T$-algebraic structure of an internal groupoid  $G$ in the semi-abelian category $\C$ lifts  to  a covering groupoid.

\begin{Thm}\label{groupoidgrouplift} Let $T$ be a semi-abelian theory with unique constant $e$,  $H$ a groupoid and $G$ an internal groupoid in $\C$ whose underlying groupoid is transitive.  Suppose that  $p\colon H\rightarrow G$ is a covering morphism of underlying groupoids, $e'\in \Ob(H)$ such that $p(e')=e$ and the characteristic group $C$ of $p$ at $e'$ is a sub $\T$-algebra of $G(e)$. Then the $\T$-algebraic structure of $G$ lifts to $H$.
\end{Thm}
\begin{Prf} Let $C$ be the characteristic group of the  covering morphism $p\colon H\rightarrow G$ at $e'\in \Ob(H)$.  By Theorem  \ref{Theoactiongpdcover}, we have a covering  morphism of groupoids $q\colon H_{C}\rightarrow G$ with the  characteristic group $C$.  Since the covering morphisms $p$ and $q$ are equivalent we can replace $p$ with $q$ and prove  that the $\T$-algebraic structure of $G$ lifts to $H_C= G\ltimes A_C$.   By Theorem \ref{semiabelactiongpd}, $A_C$ becomes a $\T$-algebra and $G$ acts on $A_C$; and by  Theorem \ref{Teointernalgpd}  $\T$-algebraic structure of  $G$  lifts to $H_C$ which completes the proof. \end{Prf}

\section*{Acknowledgement} We are grateful to  the referee for reading  the paper carefully;  and making  many constructive  comments guiding us to improve it. We would like to thank to Prof Francis Borceux for his answers to our questions via e-mail and providing the proof of Lemma \ref{Lemmapsto-e}. Our thanks is also to Prof. Ronald Brown for his useful comments.

\end{document}